\documentclass[12pt]{article}
\usepackage{amsmath,amsfonts,amssymb,amsthm}

\theoremstyle{plain}
\newtheorem{theorem}{Theorem}

\newtheorem{proposition}[theorem]{Proposition}

\theoremstyle{definition}

\newtheorem{remark}[theorem]{Remark}

\newcommand{\PP}{\mathbb{P}} 
\newcommand{\EE}{\mathbb{E}} 

\begin{document}  

\title{Quantitative version of a  Silverstein's result}
\author{ Alexander~E.~Litvak \qquad Susanna Spektor}
\date{}

\maketitle

\begin{abstract}
We prove a quantitative version of a Silverstein's Theorem on a condition for convergence in probability of the norm of random matrix. More precisely, we show that for a random matrix whose entries are i.i.d. random variables, $w_{i,j}$, satisfying certain natural conditions,
is not small with large probability.
\end{abstract}

\setcounter{page}{1}

Let $\{w_{i j}\}$  be i.i.d. random variables, identical copies of a certain random variable $w$.
Conditions on $w$ (e.g. moments, tails) will be mentioned later.
For each positive integer $n$ we consider a $p \times n$ matrix
$W_n=(w_{i j})$, $i=1, 2..., p$, $j=1, 2..., n$. We study $\lambda_{max}(\Gamma_n)$, the largest eigenvalue of
the sample covariance matrix  $\Gamma_n :=\displaystyle{\frac{1}{n} W_{n} W_{n}^T}$, of $n$ samples
of a $p$-dimensional vector containing i.i.d. components, where $p$ and $n$ are large. In the present work we considered the case when $p\leq n$, otherwise everything will work with conjugate matrices.

\medskip

It was proved in \cite{Silverstein1989}, see also \cite{Yin1988}, that with $p=p(n)$ satisfies  $\displaystyle{\frac{p}{n} \longrightarrow \beta>0 \hspace{4mm} \text{as \hspace{5mm }} n\to \infty}$ and if $\EE(w) = 0$ and $n^{4} \PP(|w| \geq n) = o(1)$, then $\lambda_{max}(\Gamma_n)$ converges in probability to the nonrandom quantity, $(1 + \sqrt{\beta})^{2} \EE(w^2)$.

\medskip

Recall, that the statement ``$\lambda_{n} := \lambda_{max}(\Gamma_n)$ converges in probability
to the limit $L:=(1 + \sqrt{\beta})^{2} \EE(w^2)$'' means that for each $\varepsilon > 0$ one has
$$
\lim_{n \to \infty} \PP\bigl( |\lambda_{n} - L| \geq \varepsilon \bigr) = 0.
$$
In particular this implies that for any $\delta > 0$ there is a positive integer $N$,
such that $n \geq N$ implies $\PP\bigl( \lambda_{n} \geq L + \varepsilon \bigr) < \delta$.

\medskip

In the present note we establish an estimate of the form
$\PP\bigl( \lambda_{max}(\Gamma_n) \geq K \bigr) \geq \delta$, where $K\geq 1$, $\delta$ is small and  dimensions of the matrix are large.

In the proposition below we use the following condition on the random variable $w$:
\begin{align}\label{1}
\forall t\geq 1, \,\alpha>0, \, c_0>0\qquad \PP(|w|\geq t)\geq \frac{c_0}{t^{\alpha}}.
\end{align}

\begin{proposition}
Let $\alpha \geq 2$, $c_0>0$. Denote by $X_i, 1\leq i\leq p$ column-vectors of $W_n$. Let also $\{w_{i j}\}$  be i.i.d. random variables, identical copies of a certain random variable $w$, such that $\EE w=0$, $\EE w^2=1$ and  satisfies condition \mbox{\emph{(\ref{1})}}. Then, for every $K \geq 1$,
\begin{align*}
\PP \left( \sup_{1\leq i \leq p} |X_i|\geq \sqrt{Kn} \right) \geq \min \left\{c_0p/(4n^{\frac{\alpha}{2}-1}K^{\frac{\alpha}{2}}), \hspace{2mm} \frac 12 \right\}.
\end{align*}
In particular,
\begin{align*}
\PP \left( \lambda_{max}(\Gamma_n) \geq K \right) \geq \min \left\{c_0p/(4n^{\frac{\alpha}{2}-1}K^{\frac{\alpha}{2}}), \hspace{2mm} \frac 12 \right\}.
\end{align*}
\end{proposition}

\begin{remark}
Note, that by Chebychev's inequality, $\displaystyle{\PP(|w|\geq t)\leq \frac{1}{t^2}}$.
Note also, that we use condition (\ref{1})  in the proof only once, with $t=\sqrt{Kn}$.
\end{remark}
\begin{remark}
 Note that if $\alpha < 4$ and $p\geq \displaystyle{{\frac{2}{c_0}}K^{\frac{\alpha}{2}}n^{\frac{\alpha}{2}-1}}$, then, by condition (\ref{1}), we have
 $\displaystyle{\frac n2 \PP(w^2\geq Kn)\geq \frac {nc_0}{2(Kn)^{\frac{\alpha}{2}}}=\frac{c_0}{2(K)^{\frac{\alpha}{2}}n^{\frac{\alpha}{2}-1}}\geq \frac 1p}$.
 Therefore,  one has $\displaystyle{\PP(\lambda_{max}(\Gamma_n)\geq K)\geq \frac 12}$ (see the Case 2 of the  proof below).
\end{remark}
\bigskip

\noindent

\medskip

\begin{proof}
For any $K\geq0$,
\begin{align}
\PP\biggl(  \sup_{1\leq i \leq p}|X_i| \geq \sqrt{Kn}  \biggr)
&= 1 - \PP\biggl( \sup_{1\leq i \leq p}|X_i| \geq \sqrt{Kn} \biggr)  \notag\\
&= 1 - \PP\biggl( \bigcap_{i=1}^{p} \{ |X_i| < K \} \biggr). \label{a1}
\end{align}

\medskip

But observe that for each $i$, $|X_i|$ depends on the entries of the $i$-th row of the matrix $W_{n} = (w_{ij})$.
Since the rows of $W_n$ are independent, it follows that $|X_i|$, for $i = 1,\ldots,p$,
are independent random variables. Also it follows that for each fixed $i$
we have that random vectors $\sum_{j=1}^n w_{i j}^2$ and $\sum_{j=1}^n w_{1j}^2$  have the
same distribution\footnote{We write $x \sim y$ to say that x and y have  the same distribution.} and, with $w_1, \ldots, w_n$ independent copies of $w$,
for each $i$ we have $\displaystyle{|X_i| \sim |X_1| \sim \frac{1}{n}\sum_{j=1}^n w_{j}^2}$.
Therefore, continuing from \eqref{a1}, we get
\begin{align}
\label{a2}
\PP\biggl(  \sup_{1\leq i \leq p}|X_i| \geq \sqrt{Kn}  \biggr)
&\geq 1 - \prod_{i=1}^{p} \PP\bigl(  |X_i| < K  \bigr)  \notag \\
&= 1-\biggl( \PP\biggl(\frac 1n \sum_{j=1}^n w_{j}^2 <  K \biggr) \biggr)^p.
\end{align}
Now we will estimate the right hand side of (\ref{a2}).

\medskip

Observe that
\begin{align*}
\PP\biggl( \frac{1}{n} \sum_{j=1}^n w_j^2 \geq K \biggr)
\geq \PP\biggl( \bigcup_{j=1}^n \{ w_j^2 \geq nK \} \biggr).
\end{align*}
For $j = 1,\ldots,n$, consider the events $A_{j} = \{ w_{j}^2 \geq nK \}$.
Since $w_{j}$'s are i.i.d. with the same distribution as $w$,
we have $\PP(A_{j}) = \PP(w^2 \geq nK )$ for all $j$ and that events $A_j$ are independent.
Using inclusion-exclusion principle, we have
\begin{align}\label{4}
\PP\biggl( \bigcup_{j=1}^n A_{j} \biggr)
&\geq \sum_{j=1}^n \PP(A_j)- \sum_{j\neq k} \PP\left(A_j\cap A_k\right)\notag \\
&=\sum_{j=1}^n \PP\left(w^2\geq nK\right) -
        \sum_{j\neq k}\bigl( \PP\left( w^2 \geq nK \right) \bigr)^2 \notag\\
&=n \PP\left(w^2\geq nK\right)-\frac{n^2-n}{2} \bigl( \PP\left(w^2\geq nK\right) \bigr)^2 \notag\\
&=\frac{n}{2} \PP(w^2\geq nK)(2-(n-1)\PP(w^2 \geq nK)).
\end{align}

Using Chebychev's inequality
$\displaystyle{\PP(w^2\geq nK) \leq \frac{1}{nK}}$, we observe  for  $K\geq 1$
$$
\displaystyle{2-(n-1)\PP(w^2\geq nK) \geq 1}.
$$

\medskip
Thus we obtain
\begin{align}\label{5}
\PP\left(\lambda_{max}(\Gamma_n)\geq K\right)
&\geq 1-\left(1 - \PP\left(\frac{1}{n} \sum_{j=1}^n w_j^2\geq K\right)\right)^p \notag \\
&\geq 1-\left(1 - \frac{n}{2} \PP\left(w^2\geq nK\right)\right)^p.
\end{align}



We consider two cases.
\medskip

\noindent{Case 1}:
\begin{align}\label{7}
\frac n2 \PP(w^2\geq Kn)\leq \frac 1p.
\end{align}
In this case, using that $\displaystyle{(1-x)^p\leq (1+px)^{-1}}$ on $[0,1]$,  we get
\begin{align}\label{8}
\PP(\lambda_{max}(\Gamma_n)\ge K)\geq 1-\frac{1}{\frac{np}{2}\PP(w^2\geq Kn)+1}.
\end{align}
Using (\ref{7}) our condition (\ref{1}) with $t=\sqrt{Kn}$, we get
\begin{align}\label{9}
1\geq \frac{np}{2}\PP(w^2\geq Kn)\geq \frac{np}{2}\frac{c_0}{(Kn)^{\frac{\alpha}{2}}}.
\end{align}
Thus,
\begin{align}\label{10}
\PP\biggl( \sup_{1\leq i \leq p}|X_i| \geq \sqrt{Kn} \biggr)\geq \frac{c_0 p}{4 n^{\frac{\alpha}{2}-1}K^{\frac{\alpha}{2}}}.
\end{align}

\medskip

\noindent {Case 2}:
\begin{align}\label{11}
\frac n2 \PP(w^2\geq Kn)\geq \frac 1p.
\end{align}
In this case, using (\ref{11}) and that $\displaystyle{(1-x)^p\leq 1-\frac{px}{2}}$, on $\displaystyle{\left[\frac 1p,1\right]}$, we get
\begin{align}\label{12}
\PP\biggl( \sup_{1\leq i \leq p}|X_i| \geq \sqrt{Kn} \biggr)\geq \frac{np}{4}\PP(w^2\geq Kn)\geq \frac 12.
\end{align}

Now, combining (\ref{10}) and (\ref{12}) we obtain
\begin{align}
\PP\biggl(\sup_{1\leq i \leq p}|X_i| \geq \sqrt{Kn} \biggr)\geq \min \left\{c_0 p/(4 n^{\frac{\alpha}{2}-1}K^{\frac{\alpha}{2}}), \quad \frac 12\right\}.
\end{align}

In particular part follows, since $\displaystyle{\lambda_{max}(\Gamma)=\|\Gamma_n\|=\frac 1n \|W_n\|^2\geq \frac 1n \sup_{1\leq i\leq p} |X_i|^2}$.

\end{proof}

\bigskip


\end{document}